\documentclass[12pt]{amsart}
\usepackage{amsmath, amssymb, amsthm, latexsym}
\usepackage{amssymb}
\usepackage{latexsym}
\usepackage[center]{caption}
\usepackage{tikz}
\usepackage{standalone}
\usepackage{graphicx}
\newcounter{braid}
\newcounter{strands}
\pgfkeyssetvalue{/tikz/braid height}{1cm}
\pgfkeyssetvalue{/tikz/braid width}{1cm}
\pgfkeyssetvalue{/tikz/braid start}{(0,0)}
\pgfkeyssetvalue{/tikz/braid colour}{black}
\pgfkeys{/tikz/strands/.code={\setcounter{strands}{#1}}}

\makeatletter
\def\cross{%
  \@ifnextchar^{\message{Got sup}\cross@sup}{\cross@sub}}

\def\cross@sup^#1_#2{\render@cross{#2}{#1}}

\def\cross@sub_#1{\@ifnextchar^{\cross@@sub{#1}}{\render@cross{#1}{1}}}

\def\cross@@sub#1^#2{\render@cross{#1}{#2}}

\def\render@cross#1#2{
  \def\strand{#1}
  \def\crossing{#2}
  \pgfmathsetmacro{\cross@y}{-\value{braid}*\braid@h}
  \pgfmathtruncatemacro{\nextstrand}{#1+1}
  \foreach \thread in {1,...,\value{strands}}
  {
    \pgfmathsetmacro{\strand@x}{\thread * \braid@w}
    \ifnum\thread=\strand
    \pgfmathsetmacro{\over@x}{\strand * \braid@w + .5*(1 - \crossing) * \braid@w}
    \pgfmathsetmacro{\under@x}{\strand * \braid@w + .5*(1 + \crossing) * \braid@w}
    \draw[braid] \pgfkeysvalueof{/tikz/braid start} +(\under@x pt,\cross@y pt) to[out=-90,in=90] +(\over@x pt,\cross@y pt -\braid@h);
    \draw[braid] \pgfkeysvalueof{/tikz/braid start} +(\over@x pt,\cross@y pt) to[out=-90,in=90] +(\under@x pt,\cross@y pt -\braid@h);
    \else
    \ifnum\thread=\nextstrand
    \else
     \draw[braid] \pgfkeysvalueof{/tikz/braid start} ++(\strand@x pt,\cross@y pt) -- ++(0,-\braid@h);
    \fi
   \fi
  }
  \stepcounter{braid}
}

\tikzset{braid/.style={double=\pgfkeysvalueof{/tikz/braid colour},double distance=1pt,line width=2pt,white}}

\newcommand{\braid}[2][]{%
  \begingroup
  \pgfkeys{/tikz/strands=2}
  \tikzset{#1}
  \pgfkeysgetvalue{/tikz/braid width}{\braid@w}
  \pgfkeysgetvalue{/tikz/braid height}{\braid@h}
  \setcounter{braid}{0}
  \let\sigma=\cross
  #2
  \endgroup
}
\makeatother

\input xypic

\swapnumbers
\newtheorem{theorem}[subsection]{Theorem}
\newtheorem{proposition}[subsection]{Proposition}
\newtheorem{lemma}[subsection]{Lemma}

\def\Proof{\medskip\noindent{\bf Proof: }}

\def\Z{\mathbb{Z}}

\def\C{\mathbb{C}}

\def\R{\mathbb{R}}
\def\C{\mathbb{C}}

\def\N{\mathbb{N}}

\def\md{\mathcal{D}}

\def\qed{\hfill$\square$\medskip}

\def\Zpk{\mathbb{Z}/p^{k}}
\def\Zpk1{\mathbb{Z}/p^{k-1}}

\newcommand{\rref}[1]{(\ref{#1})}

\newcommand{\beg}[2]{\begin{equation}\label{#1}#2\end{equation}}
\def\r{\rightarrow}
\def\mc{\mathcal{C}}

\def\sl2{\widetilde{SL_{2}(\Z)}}

\def\md{\mathcal{D}}

\title[Fields, stable homotopy and Khovanov homology]{Field theories, stable homotopy theory and Khovanov homology}
\author{Po Hu, Daniel Kriz, Igor Kriz}
\thanks{Po Hu was supported by NSF grant DMS 1104348. 
Igor Kriz was supported by NSF grant DMS 1102614 and by a Wolfensohn Fellowship at the Institute for Advanced Study}

\begin{document}

\maketitle

\section{Introduction}

\vspace{3mm}

The present paper has a somewhat peculiar history. Essentially all the work took place in Fall 2011 and
Winter 2012. It was a conglomerate of several mathematical projects. We put the outcome on the arXiv, but
neither author had strong
feelings about publication. Recently, however, interest in these topics was rekindled (see for
example \cite{lls}). We therefore decided to revise the manuscript, and publish it in the present proceedings
volume. 

As for the ingredients of the project, P. Hu and I. Kriz
were long interested in topological modular functors, i.e.
$1+1$-topological quantum field theories valued in finite-dimensional
$\C$-vector spaces, and developing a ``realization'' construction
which would convert such a structure into a $1+1$-topological quantum field theory
valued in $k$-modules, where $k$ denotes the connective K-theory spectrum.
D. Kriz, on the other hand, studied
Khovanov homology as a part of another project \cite{kk}.
From joint discussion, there arose a project of writing down a realization
construction and applying it to constructing a $k$-module refinement
of Khovanov homology.

In this, we ultimately succeeded, and we learned quite a bit. The realization
into $k$-modules is an intuitively compelling idea, but technical details
are tricky due to the difficulty of multiplicative infinite loop space
theory. We decided to use the machinery of Elmendorf and Mandell
\cite{em} which uses multicategories enriched in groupoids. We then
discovered that the multicategory language is quite a convenient tool
for axiomatizing modular functors as well. A multicategory has objects
and $n$-tuple morphisms
\beg{ei+}{X_1,\dots,X_n\r Y,
}
which compose in the same way as the elements of an operad. Multicategories are also
called colored operads or multi-sorted operads.
In this paper, by a {\em $\star$-category}, we mean a multicategory enriched in groupoids
where for every $n$-tuple $X_1,\dots,X_m$, there is a universal
morphism \rref{ei+} (in the $2$-category sense). We denote
$Y=X_1\star\dots\star X_n$. For detail, see Definition \ref{def21} below. By a $\star$-functor,
we mean a multifunctor which preserves this structure (although we
focused on the $2$-category context, there are, of course, similar concepts
in ordinary multicategories and multicategories enriched in topological
spaces). As examples of $\star$-categories, we have the $1+1$-oriented cobordism multicategory
$\mathcal{A}$ (and its many variants), and also a certain $\star$-category
$\mc_2$ associated with any symmetric bimonoidal category $\mc$ (at least
when its $2$-morphisms for a groupoid). By
a {\em $\mc$-valued modular functor} on a $\star$-category $\mathcal{Q}$
we then mean a $\star$-functor
$$\mathcal{Q}\r\mc_2.$$
Our realization theorem is then the following result (see Section \ref{smulti}
for precise definitions):

\begin{theorem}
\label{th11}
A $\mc$-valued modular functor gives rise, in a canonical way, to a multifunctor
\beg{eth11+}{M:B_2\mathcal{Q}\r \text{$k_\mc$-modules}
}
where $k_\mc$ is the $E_\infty$-ring spectrum associated with $\mc$,
and $B_2\mathcal{Q}$ denotes the topological multicategory obtained by
taking classifying spaces of the $1$-morphism groupoids.

Furthermore, a universal multimorphism
$$X_1,\dots, X_n \r X_1\star\dots\star X_n,$$
maps, under \rref{eth11+}, into an equivalence
\beg{eth112}{M_{X_1}\wedge_{k_\mc}\dots\wedge_{k_\mc}M_{X_n}\r
M_{X_1\amalg\dots\amalg X_n}.
}
\end{theorem}

\vspace{3mm}
\noindent
{\bf Comments:} 1. By ``in a canonical way'', we mean that we have a specific construction in mind. It is given by
the Elmendorf-Mandell construction.

2. The second statement requires some explanation.
What is relevant here is that we work in the category of symmetric spectra
where we have a symmetric monoidal structure under which $E_\infty$ ring spectra are, by definition,
precisely commutative monoids. For an $E_\infty$-ring spectrum $E$ in this
category, the multicategory of $E$-modules is then a $\star$-category
where $\star=\wedge_E$. The morphism \rref{eth112} is then the morphism
whose existence is the defining property of the $\star$-structure.

In view of this, it would be interesting to know if one could devise
a construction where the map \rref{eth112} would be an isomorphism instead of
just an equivalence,
i.e. such that our construction would be a $\star$-functor. Our construction
does not give this, because the Elmendorf-Mandell machine does not give
a $\star$-functor. We suspect that such a $\star$-functor might not exist.

\vspace{3mm}
The main application we had in mind was refining Khovanov ($sl_2$)-homology
of links in $S^3$ to a $k$-module invariant where $k$ is connective $k$-theory.
We hoped to achieve this by refining Khovanov's $1+1$-TQFT $\Lambda[x]$
into a $1+1$-modular functor valued in finite-dimensional $\C$-vector spaces
on the oriented $1+1$-cobordism category $\mathcal{A}$. This turns out to be impossible,
but we succeeded in constructing a modular functor on the $\star$-category $\mathcal{A}_{s}^{A}$
of spin $1+1$-cobordisms where the objects are antiperiodic $1$-manifolds.
(For a detailed definition of $\mathcal{A}_{s}^{A}$, see Example \ref{sscobord}. 
For an explanation why $\mathcal{A}_{s}^{A}$ is needed instead of $\mathcal{A}$, see
Section \ref{sec32}.)
It therefore came as a surprise when the spin-structure dropped out in the
end, and we were able to use this construction to define a $k$-theory lift
of Khovanov homology on links without spin structure. We then thought
that there must be a geometric guiding principle which explains this simplification.

\vspace{3mm}
Soon afterwards, the paper \cite{ls} by Lipshitz and Sarkar appeared on the arXiv.
This paper contains a construction of a {\em stable homotopy} refinement of Khovanov
homology. The paper \cite{ls} uses a different technique, namely Cohen-Jones-Segal
flow categories arising from Morse theory, but after some initial skepticism, we realized that Lipshitz and
Sarkar discovered the geometric principle we were looking for, while at the
same time rendering our $k$-theory refinement obsolete:
In our language, they realized that Khovanov's construction takes place
in the category enriched in groupoids $\mathcal{A}_K$ of {\em embedded} cobordisms
(in $S^2\times [0,1]$ - see Section \ref{smulti} for precise definitions).
They additionally observed what amounts to saying that the Khovanov TQFT refines into
a lax $2$-functor into $\mathcal{S}_2$ where $\mathcal{S}$ is the symmetric bimonoidal category
of finite sets. The $\star$-functor structure here is lost as $\mathcal{A}_K$ is not
a $\star$-category, but a $\star$-functor structure turns out to be unnecessary because the target of the
construction is the category of symmetric spectra (instead of modules over
another rigid ring spectrum), so the module structure does not have to be
discussed (although an analogue of \rref{eth112} is relevant and an equivalence
follows from more special arguments). We therefore end up with an alternate proof of
the following result, without requiring the language of Morse theory and flow categories:

\vspace{3mm}
\begin{theorem} (Lipshitz-Sarkar \cite{ls})
\label{t1}
There exists an explicitly defined $k$-module symmetric
spectrum $M(L)$
assigned to an oriented link $L$
such that for isotopic oriented links $L\cong L^\prime$, there exists
an equivalence
$$M(L)\simeq M(L^\prime)$$
and such that
$$M(L)\wedge H\Z$$
corresponds to the Khovanov chain complex under the equivalence of derived categories
of strict $H\Z$-modules and chain complexes \cite{ekmm}, where $S\r H\Z$ is the 
strictly commutative strict symmetric ring spectrum unit.
In other words, the homology of $M(L)$ is the Khovanov homology of $L$. In \cite{ls},
$M(L)$ is denoted by $\chi_{Kh}(L)$.
\end{theorem}

\vspace{3mm}
The convention in \cite{khovanov} is that the Khovanov complex is written as a cochain complex.
However, in our treatment, we reverse this by reversing the conventions for the $0$-resolution
and $1$-resolution of link crossings (see Figure 1 below). This has the effect of changing cohomology
into homology, which is more natural from our point of view.

\vspace{3mm}
As already remarked, strictly speaking, the full strength of Theorem \ref{th11} is
unnecessary for our proof of Theorem \ref{t1}. However, our
investigation of stable homotopy realization of modular functors, including
the construction of the Khovanov topological modular functor on $\mathcal{A}_{s}^{A}$,
provides an excellent motivation for understanding our proof of Theorem \ref{t1},
and thereby makes the argument easier to understand.
Because of this, we decided to report on both investigations in the same paper,
and also to include a discussion of the spin-dependent modular functor.

\vspace{3mm}
The present paper is organized as follows: In Section \ref{smulti},
we review the main point of the Elmendorf-Mandell formalism, and introduce
the notion of a $\star$-category and $\star$-functor, and
also prove Theorem \ref{th11}. In Section \ref{sl},
we construct our main example of the spin modular functor refinement of the Khovanov
$1+1$-dimensional TQFT $\Lambda(x)$, and also the reinterpretation
of Lipshitz-Sarkar's construction in terms of $2$-functors.
In Section \ref{scube}, we construct the refinements of the Khovanov functor. In Section \ref{sinv}, we construct 
the refined
invariant, and state a more specific version of
Theorem \ref{t1} (Theorem \ref{tinv}). Section \ref{sproof} is dedicated
to proving link invariance (Theorems \ref{t1} and \ref{tinv}), refining, essentially, the proof
of link invariance of Khovanov homology \cite{khovanov} (see also Bar-Natan \cite{dbar}).

\vspace{3mm}
\noindent
{\bf Acknowledgements:} We are indebted to Tony Licata, Robert Lipshitz and
Chris Schommer-Pries for valuable discussions.

\vspace{3mm}

\section{Multicategories and topological field theories}\label{smulti}

\vspace{3mm}

Following \cite{em}, a {\em multicategory} $\mc$ has a class of {\em objects}
$Obj(\mc)$ and classes of morphisms $Mor_n(\mc)$, $n=0,1,2,\dots$
written as
$$\phi:(x_1,\dots,x_n)\r y,\;\;\; x_1,\dots,x_n,y\in Mor(\mc).$$
We also write
$$\phi\in \mc(x_1,\dots,x_n;y).$$
There are composition, equivariance and unit axioms analogous to the definition of an operad.
Details can be found in \cite{em}. In this paper, we will be dealing with multicategories
enriched in groupoids. This means that while $Obj(\mc)$ is a class, $\mc(x_1,\dots,x_n;y)$
are groupoids, and compositions and units are functors. Associativity, unitality and
equivariance are satisfied up to natural isomorphisms, which in turn satisfy coherence
axioms modeled on cocycle conditions. Details of this context are also
fully discussed in \cite{em}. 

Therefore, we are in a $2$-categorical context. The objects of a morphism groupoid
will sometimes be referred to as $1$-morphisms, and morphisms of a
morphism groupoid as $2$-morphisms. This is the standard language of $2$-category theory.
The reader should realize that a $2$-category where the $2$-morphisms are isomorphisms is
the same thing as a category enriched in groupoids.

\vspace{3mm}
The most fundamental examples discussed in \cite{em} are the multicategory (enriched in groupoids)
$Perm$ of (small) permutative categories and the multicategory $Sym$ (enriched in topological
spaces) of symmetric spectra. In the multicategory $Sym$, morphisms $X_1,\dots,X_n\r Y$
are the same thing as morphisms
$$X_1\wedge \dots\wedge X_n\r Y$$
where $\wedge$ is the commutative, associative and unital smash product of symmetric spectra.

In some sense, the main result of \cite{em} is constructing a continuous multifunctor
$$B_2 Perm \r Sym$$
where $B_2$ means taking the classifying spaces of the $1$-morphism groupoids, thereby 
obtaining a topological multicategory. For a multicategory $M$ enriched in groupoids, let
$Sym^{M}$ denote the category of multifunctors $M\r Sym$. The other main result of \cite{em}
is their Theorem 1.4, stating that for $M$, $M^\prime$ multicategories enriched in
groupoids, and $f:M\r M^\prime$ a multifunctor that is a weak equivalence, then
the induced functor
$Sym^{M^\prime}\r Sym^M$ is a Quillen equivalence. In other words, the construction of
\cite{em} preserves weak equivalences of multicategories.

\vspace{3mm}

\subsection{Definition:}\label{def21}
A {\em $\star$-category} is a multicategory enriched in groupoids such that
for every $x_1,\dots,x_n\in Obj(\mc)$, ($n\geq 0$), there exists an
object $x_1\star\dots\star x_n$ and a $1$-morphism
$$\phi: (x_1,\dots,x_n)\r x_1\star\dots\star x_n$$
(in the case of $n=0$, one denotes the right hand side as $1$), such that for every
$1$-morphism 
$$\psi:(x_1,\dots,x_n)\r y,$$
there exists a $1$-morphism 
$$h:x_1\star\dots\star x_n\r y$$
and a $2$-isomorphism 
$$\diagram \iota:\psi\rto^\cong & h\circ\phi\enddiagram$$
and furthermore, for other such data
$$h^\prime:x_1\star\dots\star x_n\r y,$$
$$\diagram \iota^\prime:\psi\rto^\cong & h\circ\phi,\enddiagram$$
there exists a {\em unique} $2$-isomorphism 
$$\diagram \lambda:h\rto^\cong &h^\prime\enddiagram$$
such that
$$\lambda\circ Id_\phi =\iota^\prime\circ\iota^{-1}.$$
Note that for two objects $u$, $v$ satisfying the definition of $x_1\star\dots\star x_n$,
there exist $1$-morphisms $u\r v $ and $v\r u$ (unique up to $2$-isomorphism) 
whose compositions are $2$-isomorphic to the identity. 

\vspace{3mm}

In the context of multicategories enriched in groupoids, one has a notion of lax multifunctors,
analogous to lax functors of $2$-categories, where the composition and identity axioms are satisfied
up to $2$-isomorphisms satisfying the standard coherence diagrams.

\subsection{Definition:}

A {\em $\star$-functor} is a lax multifunctor $F:\mc\r\md$ between multicategories
enriched in groupoids such that $F(x_1\star\dots\star x_n)$ is a choice for
$F(x_1)\star\dots\star F(x_n)$.

\noindent
{\bf Comment:} There is a canonical $\star$-category which comes from a (lax) symmetric monoidal category:
If the symmetric monoidal structure is $\otimes$, then morphisms 
$$x_1,\dots,x_n\r y$$
are, by definition, the morphisms
$$x_1\otimes\dots\otimes x_n\r y.$$ 
This is always a $\star$-category, with 
$$x_1\star\dots\star x_n=x_1\otimes\dots\otimes x_n.$$
Not every $\star$-category, however, comes from a symmetric monoidal category.
As an example, consider the operad $A$ where $A(k)$ is the commutative monoid of non-negative
integers $(\N_0,+)$, and composition
$$A(k)\times A(n_1)\times\dots A(n_k)\r A(n_1+\dots+n_k)$$
is 
$$(x,y_1,\dots,y_k)\mapsto x+y_1+\dots+y_k-k+1.$$
The only $2$-isomorphisms are, by definition, identities.
The reader should check that this operad (and hence multicategory) satisfies
the $\star$-category axioms, but does not come from a symmetric monoidal category.

Most of the $\star$-categories discussed in this paper however come, in fact, from
(lax) symmetric monoidal categories. The reason we prefer the $\star$-category language
is that the conditions on both $\star$-categories and $\star$-functors are much simpler
to verify in comparison with symmetric monoidal $2$-categories and $2$-functors, 
since there is only a universal property to check.

\vspace{3mm}

\vspace{3mm}
\noindent
\subsection{Examples of cobordism categories:} \label{sscobord} 
1. The `basic' cobordism category $\mathcal{A}$: The objects of $\mathcal{A}$ are oriented compact smooth
oriented $1$-manifolds.
$1$-morphisms 
$$(X_1,\dots,X_n)\r Y$$ 
are oriented cobordisms between $X_1\amalg\dots\amalg X_n$
and $Y$. $2$-morphisms are orientation preserving diffeomorphisms which are identity on the boundary. 
The $\star$-category structure is given by
$$X_1\star\dots\star X_n=X_1\amalg\dots\amalg X_n$$
with the universal $1$-morphism $(X_1,\dots,X_n)\r X_1\star\dots\star X_n$ being the
identity. 
The unit object is $\emptyset$.

2. There are a number of variants of $\mathcal{A}$. One example of interest is $\mathcal{A}_s$
where $Obj(\mathcal{A}_s)$ is the class of oriented $1$-manifolds with spin structure and 
$1$-morphisms are oriented spin cobordisms between  $X_1\amalg\dots\amalg X_n$
and $Y$. Recall that a spin structure on a $1$-manifold $M$ with tangent
bundle $\tau_M$ can be
specified by giving a real bundle $\tau^{1/2}$ and
an isomorphism of real bundles
$$\tau^{1/2}\otimes_{\R}\tau^{1/2}\cong \tau.$$
An oriented circle has
two spin structures called {\em periodic} and {\em antiperiodic},
depending on whether $\tau^{1/2}$ 
is trivial or a M\"{o}bius strip. The antiperiodic spin structure
is cobordant to $\emptyset$, while the periodic one is not.
$2$-morphisms are orientation preserving diffeomorphisms which are $Id$ on the boundary including
spin, which means also identity on $\tau^{1/2}$. 
One is also interested in the $\star$-category $\mathcal{A}_{s}^{A}$ which is defined in the
same way, but one restricts to objects which are spin $1$-manifolds with antiperiodic spin
structure on each connected component. 

3. Another variant of $\mathcal{A}$ is $\mathcal{A}_K$, the embedded $1+1$- bordism category.
Objects are smooth compact $1$-dimensional submanifolds of $S^2$. $1$-morphisms
$X_1\r X_2$ are compact $2$-submanifolds $Y$ of $S^2\times [0,1]$ whose boundary is
in $S^2\times\{0,1\}$ (which $Y$ meets transversally), and such that
$Y\cap S^2\times\{0\}=X_1$, $Y\cap S^2\times\{1\}=X_2$. $2$-isomorphisms
$Y\r Y^\prime$ are isotopy classes of diffeomorphisms $\phi:S^2\times[0,1]\r
S^2\times[0,1]$ which are the identity on the boundary and restrict to diffeomorphisms
$\phi|_{Y}:Y\r Y^\prime$ (the isotopies are required to restrict to isotopies of diffeomorphisms
$Y\r Y^\prime$). Note, however, that this $2$-category has no canonical multicategory
structure. 

\vspace{3mm}
\noindent
\subsection{An example of a target $\star$-category:} 
\label{target}
Let $\mc$ be a symmetric bimonoidal groupoid. The examples we are thinking of are:
$$\begin{array}{llll}
\mc & = & R, & \parbox[t]{2in}{a commutative semiring $R$ (considered as a discrete category, i.e. the
only morphisms are identities), $+$, $\cdot$,}\\[10ex]
\mc & = & \mathcal{V}, & \parbox[t]{2in}{the category of finite-dimensional $\C$-vector spaces and isomorphisms,
$\oplus$, $\otimes$,}\\[8ex]
\mc & = & \mathcal{S}, & \parbox[t]{2in}{the category of finite sets, $\amalg$, $\times$.}
\end{array}
$$
The $\star$-category $\mc_2$ has as objects the class of all finite sets. A
$1$-morphism $(S_1,\dots,S_n)\r T$ is a $T\times (S_1\times\dots S_n)$-matrix
(thinking of $T$ as the set of rows and $(S_1\times\dots S_n)$ as the set of
columns) of objects of $\mc$. Composition is ``matrix multiplication'' with respect to
the additive and multiplicative operation of $\mc$. $2$-isomorphisms are matrices of
$\mc$-isomorphisms.

\vspace{3mm}
\subsection{Definition:} Let $\mathcal{Q}$ be a $\star$-category and let $\mc$ be a
symmetric bimonoidal category. Then a $\star$-functor $\mathcal{Q}\r \mc_2$ is called
an {\em $\mc$-valued modular functor} on $\mathcal{Q}$. We are typically interested
in examples such as $\mathcal{Q}=\mathcal{A}$, $\mathcal{Q}=\mathcal{A}_{s}^{A}$,
etc.

\vspace{3mm}
\noindent
{\bf Comment:} The idea of modular functors originates with Segal \cite{scft}, but the definition
given in \cite{scft} was not rigirous (the coherence isomorphisms were treated as equalities, thereby
neglecting the question of coherence diagrams). Hu, Kriz and Fiore developed a formalism defining
modular functors rigorously, but the formalism is awkward from the point of view of infinite loop space theory.
Hence the variant introduced in the present paper.

\vspace{3mm}
\subsection{Remark:} A $\star$-functor 
$$\mathcal{A}^s_A\r\mathcal{V}_2$$
is the flavor of $2$-vector-space valued $1+1$-dimensional
topological field theory with Spin-structure we will use in this paper. Generally speaking, 
one tends to call 
$2$-vector-space valued $1+1$-dimensional
topological field theories ``topological modular functors''. When using
that term, however, one usually considers a larger source $\star$-category
than $\mathcal{A}$. In one variant, one removes the spin structure; if there
is spin structure, one usually removes the restriction on the spin structure
on objects being antiperiodic. In the present paper, however, we are unable
to work with these notions, as the relevant examples either don't exist, or
we are unable to construct them; when constructing the $K$-theory version of
Khovanov homology, the source $\star$-category $\mathcal{A}^s_A$
is precisely what we need.

\vspace{3mm}
\begin{lemma}
\label{l1}
Let $\mc$ be a symmetric bimonoidal category. Then
there is a canonical lax multifunctor
$$\Phi:\mc_2\r Perm$$
where $Perm$ is the lax multicategory of permutative categories (see \cite{em}, Theorem 1.1 for $Perm$).
\end{lemma} 

\Proof
On objects, we set
$$X\mapsto \prod_{X} \mc.$$
On $1$-morphisms, a morphism $f:X_1\times\dots\times X_n\r Y$
in $\mathcal{C}_2$ is a $Y\times (X_1\times\dots\times X_n)$-matrix
whose entries are objects in $\mathcal{C}$. For
each $y\times(x_1,\dots,x_n)\in Y\times X_1\times\dots\times X_n$,
denote the corresponding entry by $M_{(y,x_1,\dots,x_n)}\in Obj(\mathcal{C})$. Then
$$\Phi(f):\prod_{X_1\times\dots\times X_n} \r \prod_Y\mathcal{C}$$
is given by matrix multiplication, using the $\otimes$ in $\mathcal{C}$ as the multiplication
of entries.

On $2$-morphisms, a $2$-morphism $f\Rightarrow f^\prime$ in $\mathcal{C}_2$ is a matrix
of isomorphisms $M_{(y,x_1,\dots,x_n)}\r M^\prime_{(y,x_1,\dots,x_n)}$
where $M_{(y,x_1,\dots,x_n)}$, $M^\prime_{(y,x_1,\dots,x_n)}$ are the $(y,x_1,\dots,x_n)$'th
entries of $f$, $f^\prime$, respectively. The $2$-isomorphism $\Phi(f)\Rightarrow \Phi(f^\prime)$
is induced by these isomorphisms.
\qed

\vspace{3mm}

\subsection{Construction:} \label{rmcb}
Let $\mathcal{Q}$ be a category enriched in groupoids. Denote
by $\mathcal{Q}_B$ the multicategory enriched in groupoids with objects $\{B\}\amalg Obj(\mathcal{Q})$,
where
$$\mathcal{Q}_B(\underbrace{B,\dots,B}_n;B)=E\Sigma_n,$$
(recall from \cite{em} that $E\Sigma_n$ means the torsor over $\Sigma_n$)
and for $x,y\in Obj(\mathcal{Q})$,
$$\mathcal{Q}_B(\underbrace{B,\dots,B}_m,x,\underbrace{B,\dots,B}_n;y) = E\Sigma_{m+n}\times
\mathcal{Q}(x,y).
$$
Unspecified morphism sets are empty, and composition rules are the obvious ones. 
Then the machine of \cite{em} converts a lax multifunctor
$$F:\mathcal{Q}_B\r Perm$$
into an $E_\infty$-symmetric spectrum $R$ (obtained from $F(B)$)
and an $A_\infty$-functor (associative functor in \cite{em})
$$B_2\mathcal{Q}\r E_\infty-R-\text{modules}.$$
Recall that $B_2\mathcal{Q}$ for a category $\mathcal{Q}$ enriched in groupoids is the topological category
obtained by taking the classifying space on $2$-morphisms. Theorem 1.4 of 
\cite{em} further enables us to make this {\em strict}, i.e. we obtain
a strict commutative ring symmetric spectrum $R$ and a strict functor 
$$B_2\mathcal{Q}\r R-\text{modules}.$$
Specifically, by Theorem 1.4 of \cite{em}, and $E_\infty$-ring $R$ in symmetric spectra is naturally equivalent 
to a strictly commutative ring, and an $E_\infty$-module over $R$ is naturally equivalent to a strict $R$-module.

\vspace{3mm}

\subsection{Construction:} \label{fmodules} Now let 
$$F:\mathcal{Q}\r\md$$
be a lax functor of categories enriched in groupoids, and
let $\md$ be a $\star$-category. Note that then we obtain
a canonical lax multifunctor 
$$F_B:\mathcal{Q}_B\r \md$$
given on objects by
$$F_B(x)=F(x) \;\text{for $x\in Obj(\mathcal{Q})$},$$
$$F_B(B)=1.$$
The values of $F_B$ on $1$-morphisms and $2$-morphisms are determined by
universality. If we have, in addition, a multifunctor
$$\Phi:\md\r Perm,$$
then by Construction \ref{rmcb}, we obtain a srictly commutative symmetric
ring spectrum $R$ and a strict functor
$$B_2\mathcal{Q}\r R-\text{modules}.$$

\vspace{3mm}

\subsection{Proof of Theorem \ref{th11}:}
Similar to the situation of Construction \ref{rmcb} and \ref{fmodules}, but
with extra structure: 

$\mathcal{Q}$ and
$\mc_2$ (which plays the role of $\mathcal{D}$ are $\star$-categories,
and $F$ is a $\star$-functor. 
Accordingly, we replace $\mc_B$ by a construction which takes
into account the multiplication: Let us write, say,
$$\mathcal{Q}^{alg}_{B}(\underbrace{B,\dots,B}_n;B)=E\Sigma_n,$$
$$\mathcal{Q}^{alg}_{B}(\underbrace{B,\dots,B}_n,x_1,\dots,x_n,
\underbrace{B,\dots,B}_n;y)=E\Sigma_{m+n}\times \mathcal{Q}(x_1,\dots,x_n;y).$$
Reproducing Construction \ref{fmodules} verbatim in this context,
we obtain a multifunctor
$$B_2\mathcal{Q}\r R-\text{modules},$$
as claimed (here $R=k_\mc$). This is the first statement of Theorem \ref{th11}. 

To prove the second statement, recall that while the Elmendorf-Mandell
machine does not preserve $\star$-structure, we may compose the multifunctor into
$R$-modules with a functorial cofibrant resolution, in which case it turns universal
multiplications into equivalences
$$M_{X_1}\wedge_R\dots \wedge_R M_{X_n}\r M_{X_1\star\dots\star X_n},$$
as claimed.
\qed

\vspace{5mm}

\section{A special example: Refinements of the Khovanov $\star$-functor $\mathcal{L}$}

\label{sl}

\vspace{3mm}

\subsection{Khovanov's original functor} Let us start with the ``classical'' example, i.e.
with our interpretation of Khovanov's original construction \cite{khovanov,dbar}.
Assume that $A$ is a commutative Frobenius algebra over a commutative ring $R$, i.e.
that there is an augmentation $R$-module homomorphism
$$\epsilon:A\r R$$
such that the pairing
\beg{ekhov1}{\diagram
A\otimes_R A\rto^(.6){\text{prod}} &
A\rto^\epsilon & R
\enddiagram
}
is a non-degenerate bilinear pairing over $R$. It is well known that such an $A$
gives rise to a $1+1$-dimensional TQFT where the field operators corresponding
to pairs of pants with two inbound and one outbound
(resp. two outbound and one inbound) boundary component are given
by the product and coproduct, respectively. Here the coproduct is the dual of the
product with respect to the pairing \rref{ekhov1}. In our language, at least when
$A$ is a free $R$-module on a given basis $\Lambda$,
this specifies a multifunctor
\beg{ekhov2}{\mathcal{L}:\mathcal{A}\r R_2.
}
The basis $\Lambda$ becomes the value of the multifunctor on the object $S^1$.
The example interesting from the point of view of ($SL_2$) Khovanov
homology is $A=\Lambda_\Z[x]$. In this case, let $\Lambda=\{1,x\}$ (so $x^2=0$).

The reader should be reminded that in Khovanov's construction \cite{khovanov, dbar},
the special structure of $A=\Lambda_\Z[x]$ plays a crucial role. Essentially, one needs
the sequence
$$
\diagram
0\rto & R\rto^1 &A\rto^\epsilon  & R\rto &0
\enddiagram
$$
to be exact. This is a property of $A=\Lambda_\Z[x]$ which does not happen often
in Frobenius algebras.
While analogues of Khovanov homology for other
Frobenius algebras have since been discovered (\cite{kh3,cat}), the construction
is much more involved than a straightforward analogue of the original construction
\cite{khovanov,dbar}. 

\vspace{3mm}
\subsection{Some remarks on refining the Khovanov functor to a $\mathcal{V}_2$-valued
$\star$-functor: why spin is needed}\label{sec32}

We originally tried to refine the Khovanov $\star$-functor \rref{ekhov2} to a $\star$-functor
from $\mathcal{A}$ to $\mathcal{V}_2$. We quickly realized, however, that this cannot work:
We cannot
construct a topological modular functor in the sense encountered, say, in the context of rational
conformal field theory \cite{bk,hk,fhk,scft}. One point is that in that setting, 
$\Lambda(x)$ would be the Verlinde ring of the modular functor $\mathcal{L}$.
This is generally not allowed, as the Verlinde conjecture \cite{verlinde}
asserts that the Verlinde ring, when tensored with $\C$ (i.e. the Verlinde
algebra), be semisimple, which is certainly not the case of $\Lambda(x)$. 
This is, however, not a definitive argument: while there are proofs of the Verlinde
conjecture (\cite{ms,huang,huang1}), these depend on concrete axiom systems for RCFT, which
build in semisimplicity by requiring unitarity, so a generalization 
suitable for our purposes could still exist. 

On the other hand, one can see more directly why a topological modular functor
$\mathcal{L}$ in the naive sense cannot exist: the mapping class group of a
genus $1$ oriented surface is $SL_2(\Z)$, and is generated by Dehn
twists. However, Dehn twists are required to map to trivial $2$-isomorphisms by the
$\mathcal{L}$-functor because they can be realized on an annulus $1$-morphism,
on which the value of $\mathcal{L}$ is isomorphic to the value of $\mathcal{L}$
on a unit disk, which has a trivial mapping class group (since one can always attach
a cap to one end of the annulus). On the other hand, considering
the gluing in $\mathcal{L}$ corresponding to the coproduct in $\Lambda(x)$ 
followed by the product, which gives 
\beg{eT}{1\mapsto 1\otimes x + x\otimes 1\mapsto 2x.}
Consider the non-trivial central element
$$z=\left(\begin{array}{cc}-1 &0\\0& -1\end{array}
\right)\in SL_2(\Z),$$
which corresponds to switching the two components in the middle of the gluing.
Hence, the value of $\mathcal{L}$ on $z$ must
switch the two summands corresponding to $1\otimes x$ and $x\otimes 1$ in \rref{eT},
and hence cannot be trivial, which is a contradiction. 
\vspace{3mm}
One clue was that it might actually help to replace $\mathcal{A}$ 
by $\mathcal{A}_{s}^{A}$. In the context of RCFT (\cite{ms}), modular functors are generally not
topological, as they carry an invariant called {\em central charge}. Depending
on the value of the central charge, however, the modular functors one encounters
can sometimes be made topological, depending on the
value of the central charge, by the following maneuver: one could
tensor with the inverse of modular 
functors which are {\em invertible} with respect to the tensor product (i.e.
$1$-dimensional). What invertible modular functors one encounters depends on
the exact axiomatization; a classification is given in \cite{spin}. Without
adding any structure, the invertible modular functor of
the smallest positive central charge is $Det^{\otimes 2}$, of central charge
$4$. Therefore, a modular functor of central charge
divisible by $4$ can be made topological by tensoring with a power of
$Det^{\otimes 2}$. One has
$Det$, of central charge $2$, if one allows {\em super}-structure,
i.e. $\Z/2$-grading of the modular functor. Super-structure would not be fatal
to our application, as the $\Z/2$ corresponding to the grading is known to twist $K$-theory
(cf. \cite{as}). In other words, one can replace the target 
category $\mathcal{V}$ by the category of super vector spaces (see also \cite{klai}).

\vspace{3mm}
However, even using $Det$, we can only rectify modular functors of even integral
central charge into topological ones. One can do better if one introduces spin:
There is an invertible super-modular functor of central charge $1$ which corresponds
to the $2$-dimensional chiral fermion RCFT. There is {\em not} an invertible
super-modular functor of central charge $1/2$ which would correspond to
the $1$-dimensional chiral fermion, but a part of the modular functor
restricted to $\mathcal{A}_{s}^{A}$ (i.e. boundary components with anti-periodic
spin structure) does exist (cf. \cite{spin}), and moreover, on this
restriction to $\mathcal{A}_{s}^{A}$,
the super-structure trivializes.

\vspace{3mm}
Of course, since we have not constructed an RCFT in any generalized sense which
would correspond to $\mathcal{L}$, so we do not know what its central charge
would be. However, we see that spin can help in making the functor topological,
as long as the central charge is a multiple of $1/2$, and as long as we restrict
to $\mathcal{A}$. We do not know if the restriction to $\mathcal{A}_{s}^{A}$ is
necessary when defining a $\mathcal{V}_2$-refinement of $\mathcal{L}$, as constructing a modular functor
with spin including periodic boundary components is much harder to
do ``by hand''.

\vspace{3mm}

\subsection{A $\mathcal{V}_2$-refinement of the Khovanov $\star$-functor} 
We will now construct ``by hand'' a certain lax $\star$-functor
\beg{ekhov3}{\mathcal{L}_s:\mathcal{A}^s_A\r\mathcal{V}_2.}
On objects, let $C$ be a closed $1$-manifold with spin structure such that
every connected component is anti-periodic. Denote the set of connected components
of $C$ by $\pi_0(C)$. Then let
\beg{ekhov3a}{\mathcal{L}_s(C)=\prod_{c\in\pi_0(C)} \{1,x\}.}
Before specifying the effect of $\mathcal{L}$ on $1$-morphisms and $2$-morphisms,
it is helpful to introduce the following terminology for boundary components of
a compact oriented surface $\Sigma$ with spin structure, whose boundary components
are labelled $1$ or $x$: A {\em true inbound} boundary component is an inbound
boundary component labelled $1$ or an outbound boundary component labelled $x$.
A {\em true outbound} boundary component is an outbound boundary component labelled
$1$ or an inbound boundary component labelled $x$. 

\vspace{3mm}

Now for a ($2$-dimensional) oriented spin cobordism $\Sigma$ with antiperiodic
boundary components, define $\mathcal{L}_s(\Sigma)$ as follows. Let $to(\Sigma)$
denote the number of true outbound boundary components of $\Sigma$, and let
$g(\Sigma)$ denote the genus of $\Sigma$.

If $\Sigma$ is connected, then
\beg{ekhov+}{\mathcal{L}_s(\Sigma)=\left\{\begin{array}{ll}
\C &\text{if $g(\Sigma)=0$ and $to(\Sigma)=1$}\\
\C\oplus\C &\text{if $g(\Sigma)=1$ and $to(\Sigma)=0$}\\
0 & \text{else.}
\end{array}\right.}
By definition of a $\star$-functor, we must, of course, for a general cobordism $\Sigma$, have
$$\mathcal{L}_s(\Sigma)=\bigotimes_{\Sigma^\prime}\mathcal{L}(\Sigma^\prime)$$
where $\Sigma^\prime$ runs through the connected components of $\Sigma$.

\vspace{3mm}
\begin{lemma}
This defines a lax $\star$-functor
$$\mathcal{L}_s:\mathcal{A}^s_A\r \mathcal{V}_2.$$
\end{lemma}

\begin{proof}
Check the axioms.
\end{proof}

\vspace{3mm}	
Regarding $2$-isomorphisms, any $2$-isomorphism between spin cobordisms 
of genus $0$ is sent to the identity. To go further, it
is convenient to introduce some terminology.
By a {\em reference curve} in a genus $1$ Kervaire invariant $0$ ($2$-dimensional)
spin cobordism $\Sigma$ with antiperiodic boundary components only, we mean
an isotopy class of non-separating antiperiodic closed oriented curves in $\Sigma$.
Let $\check{\Sigma}$ denote the surface obtained from $\Sigma$ by gluing disks
to all boundary components. Without loss of generality, a reference curve 
$\alpha_\Sigma$ is chosen in each Kervaire invariant $0$ genus $1$ spin cobordism
$\Sigma$ with antiperiodic boundary components. 

Now let $f:\Sigma\r T$ be a $2$-morphism where $\Sigma$, $T$ are of genus $1$,
Kervaire invariant $0$. Let $\alpha\in H_1(\check{T},\Z)$ be the homology class
represented by $\alpha_T$. Let $(\alpha, \beta)$ be any ordered basis of
$H_1(\check{T},\Z)$ containing $\alpha$. Let 
$$f(\alpha_{\Sigma})=k\alpha+\ell\beta\in H_1(\check{T},\Z).$$
Then
\beg{ekhov++}{
\mathcal{L}_s(f)=\left\{
\begin{array}{ll}
\left(\begin{array}{rr}
1&0\\0&1\end{array}
\right) &\text{if $k\equiv 1\mod 4$, $\ell\equiv 0\mod 2$}\\[4ex]
\left(\begin{array}{rr}
0&1\\1&0\end{array}
\right) &\text{if $k\equiv -1\mod 4$, $\ell\equiv 0\mod 2$}\\[4ex]
\left(\begin{array}{rr}
\frac{1+i}{2}&\frac{1-i}{2}\\ \frac{1-i}{2}&\frac{1+i}{2}\end{array}
\right) &\text{if $k\equiv 0\mod 2$, $\ell\equiv 1\mod 4$}\\[4ex]
\left(\begin{array}{rr}
\frac{1-i}{2}&\frac{1+i}{2}\\ \frac{1+i}{2}&\frac{1-i}{2}\end{array}
\right) &\text{if $k\equiv 0\mod 2$, $\ell\equiv -1\mod 4$}
\end{array} 
\right.
}
It is easy to show that those are the only possibilities for $k$, $\ell$,
and that the definition does not depend on the choice of $\beta$.

Note that in all other (connected) cases of $f:\Sigma\r T$, $\mathcal{L}_s(f):0\r 0$,
so there is no choice in that case.

\vspace{3mm}
This does not quite conclude the definition of $\mathcal{L}_s$. Since $\mathcal{L}_s$ is
a lax multifunctor, we must specify a $2$-morphism
$$\mathcal{L}_s(f)\circ(\mathcal{L}_s(g_1),\dots,\mathcal{L}_s(g_n))\r\mathcal{L}_s(f\circ
(g_1,\dots,g_n))$$
where applicable. As it turns out,the only non-trivial case occurs when we are gluing
genus $0$ connected cobordisms $\Sigma$, $\Sigma^\prime$ into a genus $1$ connected
cobordism. In this case, let 
$$1\mapsto \left(\begin{array}{r}1\\0\end{array}\right)$$
if the true outbound boundary component $c$ of $\Sigma$ (or, equivalently, $\Sigma^\prime$)
maps (with orientation) to $\alpha\in H_1(T,\Z)$.

It then follows from the structure that if
$$c\mapsto k\alpha+\ell\beta\in H_1(T,\Z),$$
then
$$
1\mapsto
\left\{\begin{array}{ll}
\left(\begin{array}{r}
1\\0
\end{array}\right) &\text{if $k\equiv 1\mod 4$, $\ell\equiv 0\mod 2$}\\[4ex]
\left(\begin{array}{r}
0\\1
\end{array}\right) &\text{if $k\equiv -1\mod 4$, $\ell\equiv 0\mod 2$}\\[4ex]
\left(\begin{array}{r}
\frac{1+i}{2}\\\frac{1-i}{2}
\end{array}\right) &\text{if $k\equiv 0\mod 2$, $\ell\equiv 1\mod 4$}\\[4ex]
\left(\begin{array}{r}
\frac{1-i}{2}\\\frac{1+i}{2}
\end{array}\right) &\text{if $k\equiv 0\mod 2$, $\ell\equiv -1\mod 4$}\\[4ex]
\end{array}
\right.
$$
and that no other possibility can arise. 

\vspace{3mm}
\noindent
{\bf Remark:}
It is possible to use the functor $\mathcal{L}_s$ to define a $k$-module
refinement of Khovanov homology. When we did this in the original
version of this paper, however, we eventually observed that spin completely
drops out of the story (by a mechanism which we will briefly describe below).
This is the effect of a geometric principle which we will now discuss.

\vspace{3mm}

\subsection{The Lipshitz-Sarkar refinement of the Khovanov functor}

What is in fact happening is that it suffices to construct a ``field theory'' on
$\mathcal{A}_K$, i.e. an ``embedded field theory''. Indeed, reinterpreting
the construction of Lipshitz and Sarkar \cite{ls}, one can construct a lax functor
\beg{ekhov4}{\mathcal{L}_K:\mathcal{A}_K\r \mathcal{S}_2
}
(see Section \ref{target} for the definition of $\mathcal{S}_2$). 
Note again that $\mathcal{A}_K$ is not a $\star$-category
so we lose the possibility of a $\star$-structure, but on the other hand,
composing with the Elmendorf-Mandell machine (or, alternately, essentially
any infinite loop space machine which lands in symmetric spectra), we
obtain a functor
$$B_2\mathcal{A}_K\r Sym,$$
which is sufficient, since symmetric spectra are the same thing as modules over
the sphere spectrum in that category.

The construction of \rref{ekhov4} is, in a way, similar to the construction of
\rref{ekhov3}. On objects, use the same definition as for $\mathcal{L}_s$
(see \rref{ekhov3a}). On $1$-morphisms, we also adapt the definition
\rref{ekhov+}: For a connected $1+1$-cobordism $\Sigma$ embedded
in $S^2\times [0,1]$ whose boundary is in $S^2\times\{0,1\}$ which $\Sigma$
meets transversally, we set
\beg{ekhov+a}{\mathcal{L}_K(\Sigma)=\left\{\begin{array}{ll}
\{1\} &\text{if $g(\Sigma)=0$ and $to(\Sigma)=1$}\\
\{1,2\} &\text{if $g(\Sigma)=1$ and $to(\Sigma)=0$}\\
\emptyset & \text{else.}
\end{array}\right.}
In general, we set
$$\mathcal{L}_K(\Sigma)=\prod_{\Sigma^\prime}\mathcal{L}_K(\Sigma^\prime)$$
where $\Sigma^\prime$ runs through the connected components of $\Sigma$.

But how can we make consistent choices of $\mathcal{L}_K$ on $2$-morphisms
when the ``square root'' of the transposition map $c:\{1,2\}\r\{1,2\}$ cannot
be a map of sets, and only exists as a morphism of $\C$-vector spaces?

\vspace{3mm}
\noindent
{\bf Remark:}
The answer is at the heart of the problem, and was essentially discovered by
Lipshitz and Sarkar \cite{ls} in their concept of ladybug matching. In the language
of the present paper, the point is that embedding into $S^2\times [0,1]$ restricts
modular transformations severely. In fact, the embedded mapping class group
of an unknotted torus $T$ embedded in $S^2\times [0,1]$ is $\Z/2$. For, if we choose the
reference curves $\alpha, \beta$ to be fundamental cycles representing the inside
and outside of $T$, then $\alpha$ and $\beta$ must be preserved up to orientation,
and their orientations must be either both preserved or both reversed to preserve
the orientation of $T$.

\vspace{3mm}
If $\sigma$ is the generator of this $\Z/2$, we define
\beg{ekhov5}{\mathcal{L}_K(\sigma)=\tau.}
Finally, we must define the composition isomorphism when gluing two genus $0$
embedded connected cobordisms $\Sigma$, $\Sigma^\prime$ into a genus $1$
connected cobordism. In this case, let
$$1\mapsto 1$$
if the true outbound boundary component $c$ of $\Sigma$ (or, equivalently, of
$\Sigma^\prime$) maps, with orientation, to $\alpha$ or $\beta$, and
$$1\mapsto 2$$
if $c$ maps to $-\alpha$ or $-\beta$. This definition depends on the choice of
orientations of
the generators $\alpha$ and $\beta$ which indicates $4$ possible choices, but
we also have the possibility of simultaneously reversing the orientations
of $\alpha$ and $\beta$ (i.e. applying the modular transformation $\sigma$), which 
equates two and two of the choices. Therefore, there are two intrinsically different
choices to make, which corresponds to the left and right ladybug matchings of \cite{ls}.
(Also see \cite{lls} for another description.)

\vspace{3mm}

\section{The Khovanov cube functor}
\label{scube}

\subsection{Lax categories}
\label{laxcat}

Let $\mc$ be a small category. We define a category $\mc^\prime$ enriched in groupoids
where
$$Obj (\mc^\prime)=Obj(\mc),$$
$$Mor_1(\mc^\prime)=\Gamma\mc$$
where $\Gamma$ denotes the free category on a directed graph (a directed graph is a 
pair of maps $S, T$ from a set of {\em morphisms} to a set of {\em objects}). Here, we
regard $\mathcal{C}$ as a graph by forgetting that compositions exist.

There is a canonical functor 
$$\theta:\Gamma\mc\r\mc$$
(the monad structure).
There is a single $2$-isomorphism in $\mc^\prime$ between any two morphisms whose
images under $\theta$ coincide. 

\vspace{3mm}
\subsection{Links and link cobordisms}

Let $L$ be a link with spin structure, and let $\md$ be a non-degenerate projection
of $L$, i.e. an immersion into $S^2$
with only at most transverse double points (i.e. where crossings occur at angles $\neq 0,\pi$). 
Label the crossings of $\mathcal{D}$ by $1,2,\dots,n$. For the $i$'th crossing, select a disk
$D_i$ which is a neighborhood of the crossing, such that $D_1,\dots,D_n$ are disjoint. Recall that
for $\epsilon=0,1$, the $\epsilon$-resolution is obtained by replacing a chosen crossing by a non-crossing
according to Figure 1:

 \begin{figure}[!ht]
 \includegraphics{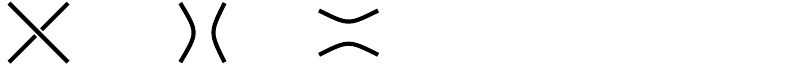}
\caption{A crossing, the $0$-resolution, 
the $1$-resolution}
\label{fskein}
\end{figure}

Recall from \cite{khovanov} that the {\em link cobordism} $\Sigma_\mathcal{D}$ is obtained by
taking 
\beg{ecube+}{\left(\md\smallsetminus\left(\bigcup_{i=1}^{n}D_i\cap \mathcal{D}\right)\right)\times I}
and for each crossing, gluing in an $\epsilon$-resolution of the crossing
at $(D_i\cap \mathcal{D})\times\{\epsilon\}$ for $\epsilon=0,1$, and
a saddle between the two crossings in $(D_i\cap\mathcal{D})\times I$.

Observe that $\Sigma_\mathcal{D}$ can be obtained by taking a ribbon
along $L$ which takes a $1/2$-twist at each crossing (thus creating
a horizontal square) and is vertical elsewhere,
and identifying the two horizontal squares over each crossing. Note that
the ribbon always has an even number of half-twists, since there are two
per crossing. Hence, the ribbon may be identified with $\tau_L\otimes_\R \C$. 

\vspace{3mm}
\noindent
{\bf Comment:} To avoid confusion, note that in the present paper, by a link cobordism
we mean the surface ($1+1$-embedded cobordism) associated with a link projection,
{\em not} a cobordism of links.

\vspace{3mm}
Now let us observe that complete resolutions of a link projection are, by definition,
objects of $\mathcal{A}_K$, and the cobordisms $\Sigma_L$ are $1$-morphisms on
$\mathcal{A}_K$. Let us also make another observation: Let $\mathcal{D}$ be
a non-degenerate link projection; label its crossings $D_1,\dots, D_n$.
Let $\mathcal{D}^\prime$  be the projection obtained by taking $0$-resolutions of
$D_1,\dots,D_k$ and let $\mathcal{D}^{\prime\prime}$ be the projection obtained
by taking $1$-resolutions of $D_{k+1},\dots,D_n$. Then there is a canonical 
$2$-isomorphism in $\mathcal{A}_K$
\beg{ecube++}{\diagram\Sigma_{\mathcal{D}^{\prime\prime}}\circ\Sigma_{\mathcal{D}^{\prime}}
\rto^\cong & \Sigma_{\mathcal{D}}\enddiagram}
We call these $2$-isomorphisms {\em gluing isomorphisms}.

\subsection{Definition of the lax cube functor}
\label{laxcube}

Now let $I$ be the category with two objects $0,1$ and a single morphism $0\r 1$
(and no morphism $1\r 0$).
Now let $\md$ be an admissible link projection with $n$ crossings.
Then $\md$ determines a lax functor
$$\mc:({I}^n)^\prime\r \mathcal{A}_K$$
(the category $(I_n)^\prime$ is defined in Section \ref{laxcat}) as follows:

To define $\mc$ on objects,
$$(\epsilon_1,\dots,\epsilon_n), \;\epsilon_i\in\{0,1\},$$
maps to the complete resolution of $\md$ obtained by taking the
$\epsilon_i$-resolution at the $i$'th crossing.

On $1$-morphisms, consider the $I^n$-morphism
$$(\epsilon_1,\dots,\epsilon_n)\r (\epsilon_{1}^{\prime},\dots,\epsilon_{n}^{\prime}).$$
Let $J\subseteq\{1,\dots,n\}$ be the subset such that $\epsilon_{i}^{\prime}=\epsilon_{i}+1$ for $i\in J$ and 
$\epsilon_{i}^{\prime}=\epsilon_{i}$ for $i\notin J$. Then this
$I^n$-morphism is sent to 
$$\Sigma_{\md^\prime}$$
where $\md^\prime$ is the projection obtained from $\md$ by taking the 
$\epsilon_i$-resolution at the $i$'th crossing for all $i\notin J$.

$2$-morphisms of $(I^n)^\prime$ are sent to gluing isomorphism
of the surfaces $\Sigma_{\md^\prime}$ for different $\md^\prime$, by their
shared boundary component.

\vspace{3mm}
\subsection{The spin data}

In this subsection, we will discuss directly spin structures on link
projections. While this gives additional geometric insight into the $k$-theory
refinement of Khovanov homology, this material is not strictly necessary to follow
the progression of the paper, and the reader only interested in the
proof of Theorem \ref{t1}  may skip it. 

By a {\em link with spin structure}, we mean a real bundle $\tau^{1/2}_{L}$
together with an isomorphism
$$\tau^{1/2}_{L}\otimes_\R\tau^{1/2}_{L}\r \tau_L$$
where $\tau_L$ is the tangent bundle of $L$. Note that this specifies an orientation
on $\tau_L$ where we call a tangent vector {\em positive} if it has a square
root in $\tau^{1/2}_{L}$.

By a {\em projection with spin} we mean a non-degenerate projection $\md$ of $L$
together with a {\em spin} of the self-identification of the ribbon
$\tau_L\otimes_\R\C$ along each crossing square, namely an automorphism of
the bundle $\tau^{1/2}_{L}\otimes_\R\C$ which covers the identity on 
$\tau_L\otimes_\R\C$. By gluing of bundles, this data, given in a projection with spin, specifies a spin
structure on $\Sigma_{\mathcal{D}}$.

\vspace{3mm}
Recall that for a complex $1$-manifold $\Sigma$ with spin (i.e. a complex line
bundle $\tau^{1/2}_{\Sigma}$ and an isomorphism $\tau^{1/2}_{\Sigma}\otimes_\C
\tau^{1/2}_{\Sigma}\cong \tau_\Sigma$), and an oriented curve $c$ in $\Sigma$,
we have a determined spin structure on $c$ where $(\tau^{1/2}_{c})_x$ is spanned
by the $(\tau^{1/2}_{\Sigma})_x$-square roots of a positive tangent vector
to $c$ at $x$. 

\vspace{3mm}
We call a projection $\md$ with spin of $L$ {\em admissible} if the induced spin
structure on every non-self-intersecting circuit in $\mathcal{D}$ is antiperiodic. (Recall \ref{sscobord}, 2.)
It suffices to verify this condition for faces. 

\vspace{3mm}

Now there is an obvious way (by sliding) to give spin-structure to R2- and
R3-moves. R1-moves require a more detailed discussion, as they do interfere with
spin. When making an R1-move, we create a new face which borders the edge created
by the R1-move only. Since we will be primarily interested in admissible projections
with spin, we will only be interested in R1-moves where the new face has an antiperiodic
spin structure. Given this condition, there are two possible ways of 
introducing a
spin structure on the projection after the R1-move: one does not change the 
spin-structure on the link $L$, but changes the spin structure on the two faces
previously adjacent to the arc on which we performed the R1-move: we will call
this an R1L-move. Taking an R1L-move and changing the spin-structure of the
resulting projection by reversing the gluing of the spin structure on the new
crossing and also in the middle of the new arc created by the move, we obtain
a move which does not change the spin structures of any of the faces of the old
projection, but reverses the spin-structure of the connected component of $L$
on which we performed the move. We will call this an R1A-move.

Note that two R1L-moves on the same arc of a projection with spin is the
same as a pair of R1A-moves on the same arc: The resulting move changes
neither the spin structures of any of the faces of the old projection, nor
the spin structure of the link. We will call such a pair {\em a pair
of adjacent R1A-moves}.

\vspace{3mm}

\begin{lemma}\label{ladmex}
An admissible projection with spin of a link $L$ with spin always exists.
\end{lemma}

\Proof
Start with any projection with spin. Making $\{A,P\}$ into a group by making $A$
the neutral element, the spin structure of the infinite face is the product of 
the spin structures of the finite faces, and hence there are an even number of 
$P$-faces, including the infinite face. This specifies a $\Z/2$-valued 
$0$-cycle $\zeta$ on the CW-decomposition of $S^2$ dual to $\md$, such that
the augmentation of $\zeta$ is $0$, and hence $\zeta$ is a boundary,
$\zeta=dc$ for some $\Z/2$-valued $1$-chain $c$. The $1$-chains of $\md$ and
its dual are the same; perform an R1L-move on each arc of $\md$ on which $c$ has
coefficient $1$.
\qed

\vspace{3mm}
In fact, we have a stronger statement:

\begin{lemma}
\label{lunique}
Consider a non-degenerate projection $\md$ (without spin) of a link $L$. Then 
an admissible spin structure on $\mathcal{D}$ always exists and
any two admissible spin structures on $\md$ (for any spin
structure on $L$) are isomorphic. In particular, the spin structure on $L$ is
determined. 
\end{lemma}

\Proof
Consider the link cobordism $\Sigma_\md$ associated
with $\md$. Then $\Sigma_\md$ is an oriented surface, so the embedding
$\Sigma_\md\subset \R^3$ extends to an embedding $\Sigma_\md\times I\r \R^3$,
and the spin structure extends, of course, uniquely to $\Sigma_\md \times I$.
Smooth out $\Sigma_\md\times I$ into a manifold with boundary $\widetilde{\Sigma_\md\times
I}$. Now since the spin structure on $\Sigma_\md$ is admissible, we may attach
a disk $D_f$ to each face $f$ of $\md$ in $\Sigma_\md$, and extend the spin
structure. Hence, we may attach a copy of $D_f\times I$ to $\widetilde{
\Sigma_\md\times I}$ (and again smooth) for each face $f$ of 
$\md$, and extend the spin structure to the resulting manifold $\Gamma$
with boundary. 

The manifold $\Gamma$, however, is diffeomorphic to $D^3$, and hence has
a unique spin structure (up to isomorphism). This means that any two admissible
spin structures on $\Sigma_\md$ are isomorphic. 

Conversely, the same construction also implies that an admissible spin structure always 
exists.
\qed

\begin{proposition}
\label{pcow}
Two admissible projections with 
spin $\md$, $\md^\prime$ represent isomorphic links with spin if and only if
they are related by R2-moves, R3-moves and pairs of adjacent R1-moves. Two admissible
projections with spin $\md$, $\md^\prime$ represent isomorphic links without spin
if and only if they are related by R2-moves, R3-moves and R1A-moves.
\end{proposition}

\Proof
Consider first the second statement. Sufficiency is obvious, as the Reidemeister
moves do not change the isomorphism class of the link (without spin).
To prove necessity, suppose $\md$, $\md^\prime$ are admissible projections which
represent isomorphic
links (without spin). As is well known, disregarding spin, $\md$ can be converted to
$\md^\prime$ by a sequence of R1-moves, R2-moves and R3-moves. Now we may
give spin to the moves (preserving admissibility) by interpreting the R1-moves
as R1A-moves. By Lemma \ref{lunique}, the admissible spin structure on 
$\md^\prime$ obtained by the moves is the same as the admissible spin structure
originally given. 

\vspace{3mm}
Now consider the first statement (on links with spin structure). Again,
sufficiency is obvious as R2-moves, R3-moves and pairs of adjacent R1-moves
do not change the spin structure of the underlying link. To prove necessity, 
suppose $\md$ and $\md^\prime$ represent the same link with spin structure.
Proceed in the same way as in the part of the statement on links without
spin. Note in particular that the argument there did not depend on
the order of the Reidemeister moves chosen. By Coward's theorem
\cite{coward}, we may choose the moves in such a way that all the R1 moves come
first, followed by R2-moves, R3-moves and reversed R2-moves. Now since, when
considering spin, we interpret the R1-moves as R1A-moves, there must be an
even number of such moves on each connected component of the link in order for
the spin structures on the links corresponding to $\md$ and $\md^\prime$ to
be the same. However, note that a pair of R1A-moves on the same connected component
of a link $L$ can always be obtained as a pair of adjacent R1-moves, followed by
R3 and R2 (and possibly reversed R2) moves.
\qed

\vspace{3mm}

Analogously with \ref{laxcube}, an admissible link projection with spin with $n$ 
crossings now directly determines a lax functor
$$\mathcal{C}_s:(I^n)^\prime\r \mathcal{A}_{s}^{A}.$$
Note that for $\mathcal{D}$, $\mathcal{D}^\prime$ as above, we have a unique, up to isotopy, inclusion
$$\Sigma_{\md^\prime}\subseteq \Sigma_\md$$
commuting with the projection to $\R^2$, so there is a canonical spin structure
on $\Sigma_{\md^\prime}$ induced from the spin structure on $\Sigma_{\md}$.

\vspace{3mm}
Note that Lemma \ref{lunique} also implies that there is a canonical lax functor
\beg{ecube*}{\mathcal{A}_K\r \mathcal{A}_{s}^{A},
}
so we could simply obtain $\mathcal{L}_s$ as a composition of
\rref{ecube*} with $\mathcal{L}_K$. This way, however, we lose the
$\star$-structure, since $\mathcal{A}_K$ is not a $\star$-category and
\rref{ecube*} is not a $\star$-functor.

\vspace{3mm}

\section{Stable homotopy realization, and link invariance}
\label{sinv}

\vspace{3mm}
\subsection{}
\label{sdelta}
Now let $\md$ be an admissible projection of a link $L$. (Note: Spin structure
is not used in this Section.)
In \ref{laxcube}, we constructed a lax functor
$$\mc: (I^n)^\prime \r \mathcal{A}_K.$$
In Section \ref{sl}, we constructed a lax functor 
$$\mathcal{L}_K:\mathcal{A_K}\r \mathcal{S}_2.$$
In Lemma \ref{l1}, we further constructed a lax multifunctor
$$\Phi:\mathcal{S}_2\r Perm.$$
By the remark at the end of Subsection \ref{rmcb}, then, the composition $\Phi\mathcal{L}_K\mathcal{C}$ is canonically 
converted into a {\em strict} functor
\beg{estable+}{\Delta_\md:I^n\r Sym}

\vspace{3mm}
\noindent
{\bf Remark:} We may of course smash the functor \rref{estable+} with $k$
in the category of symmetric spectra. Alternately, we may directly consider
the composition
\beg{estable++}{(I^n)^\prime\r\mathcal{A}_K\r \mathcal{A}_{s}^{A}.
}
By Construction \ref{fmodules}, the $\star$-functor 
$$\mathcal{L}_{s}:\mathcal{A}_{s}^{A}\r\mathcal{V}_2$$
determines a lax functor
$$\diagram (I^n)^{\prime}_{B}\rto &\mathcal{V}_2\rto^{\Phi}& Perm\enddiagram$$
which, by construction \ref{rmcb}, gives a strict functor
$$I^n\r R-\text{modules}$$
where $R$ is the strictly commutative symmetric ring spectrum arising
by the Elmendorf-Mandell machine \cite{em}
from the bipermutative category $\Phi(1)$. However, $\Phi(1)$ is the category $\mathcal{V}$
of finite-dimensional complex vector spaces and isomorphisms with its usual
bipermutative category structure, so $R$ is $k$, the connective $k$-theory 
strictly commutative symmetric ring spectrum. We have, therefore, constructed
a strict functor
\beg{edelta}{\Delta_{\md,s}:I^n\r k-modules.}
While this direct construction contributes nothing to Theorem \ref{t1} as stated,
it is interesting to note that it shows that the $k$-theory realization ``remembers less data''
about the structure of the link, since it only depends on the composition \rref{estable++},
and not the embedded link cobordism.

\vspace{3mm}
\subsection{The higher cofiber}

The higher cofiber is a functor $C_n$ from the category of diagrams
$$\Gamma:I^n\r R-modules$$
to the category of $R$-modules where $R$ is a strictly commutative 
symmetric ring spetrum. Functors of such form are used extensively, for example,
in Goodwillie calculus. (See, for example, \cite{jardine} for an overview of such functors.)

One description of the higher fiber proceeds as follows. Consider the category
$\mathcal{I}$ whose objects are functions $\phi:J\r\{0,1\}$ where $J
\subseteq \{1,\dots,n\}$ and there is a unique morphism $\phi\r\psi$ if and
only if $\phi$ is a restriction of $\psi$. In other words, $\mathcal{I}$ can be thought of
as 
$$(0\leftarrow \cdot\rightarrow 1)^n.$$
Then $\Gamma$ specifies a functor
$$\widetilde{\Gamma}:\mathcal{I}\r R-modules$$
where 
$$\widetilde{\Gamma}(\phi)=\left\{\begin{array}{ll}* & \text{if $0\in Im(\phi)$}\\
\Gamma(1-\chi_J) & \text{else}
\end{array}\right.
$$
where $\chi_J(x)=1 $ if $x\in J $ and $\chi_J(x)=0$ if $x\notin J$. 
The value of $\widetilde{\Gamma}$ on morphisms is given by the corresponding
morphism values of $\Gamma$ when the target is not $*$, and by the trivial
map else. 

One defines
$$C_n\Gamma= hocolim \widetilde{\Gamma}.$$
(The right-hand side is well defined using the simplicial realization in 
$R$-modules.)

\vspace{3mm}
The advantage of the above description is that it is obviously symmetrical in 
the coordinates. There is an alternate elementary inductive 
description which is not symmetrical
in coordinates, but symmetry is readily proved by equivalence with the above description:

We define $C_0\Gamma=\Gamma$. Assuming we have already defined $C_{n-1}$, define
$$\Gamma_\epsilon:I^{n-1}\r R-modules,\;\epsilon=0,1$$
by
$$\Gamma_\epsilon=\Gamma(?,\dots,?,\epsilon).$$
Then $\Gamma$ gives a natural transformation
$$\iota:\Gamma_0\r\Gamma_1.$$
Inductively, we get a natural transformation
$$C_{n-1}\iota:C_{n-1}\Gamma_0\r C_{n-1}\Gamma_1.$$
Let $C_n\Gamma$ be the homotopy cofiber of $C_{n-1}\iota$.

\vspace{3mm}
In fact, the symmetric description of the higher cofiber immediately gives
the following fact, which will be useful to us:

\vspace{3mm}

\begin{lemma}
\label{lcofiber}
Let $\Gamma:I^n\r R-modules$ be a functor,
and let $f:\{1,\dots,k\}\r\{1,\dots,n\}$, $g:\{1,\dots,m\}\r\{1,\dots,n\}$
be maps such that $f\amalg g$ is a bijection (so, in particular, $n=k+m$).
Define for $\phi:\{1,\dots,k\}\r\{0,1\}$, a functor $\Gamma_\phi:I^m\r R-modules$
by $\Gamma_\phi(\psi)=\Gamma_{(\phi\amalg \psi)(f\amalg g)^{-1} }$ for
$\psi:\{1,\dots,m\}\r\{0,1\}$. Then $C_m\Gamma_?: I^k\r R-modules$ is
a functor in the obvious way. We have
$$C_k(C_m\Gamma_?)=C_n\Gamma.$$
\end{lemma}
\qed

\vspace{3mm}
From now on, we shall work only with the Lipshitz-Sarkar realization, i.e.
in the category $Sym$ of symmetric spectra. Analogous results in $k$-modules
follow by applying $?\wedge k$ or alternately using analogous reasoning
directly for the $k$-module realization.

\vspace{3mm}
Now recalling \rref{edelta}, we can assign, to an admissible projection
$\md$ of a link $L$ a $k$-module $C_n\Delta_\md$.

\vspace{3mm}

\begin{theorem}
\label{tinv}
If $\md$, $\md^\prime$ are nondegenerate projections of a link $L$,
then there exists an equivalence of symmetric spectra
$$\Sigma^{-n_-(\md)}C_n\Delta_\md \simeq \Sigma^{-n_-(\md^\prime)}C_n\Delta_{\md^\prime}$$
where $n_-(\md)$ denotes the number of negative crossings of the projection $\md$
(a number which does not depend on spin).
\end{theorem}

\vspace{3mm}

\section{Proof of the main theorem}

\label{sproof}

The proof of theorem \ref{tinv} basically mimics Khovanov's proof of the
invariance of Khovanov homology (see \cite{khovanov,dbar}). Of course, we cannot refer to elements
and take chain differentials; we must phrase everything in the language
of categories. We begin with two lemmas on higher cofibers:

\vspace{3mm}

\begin{lemma}
\label{lbar1}
Consider a diagram $M$ of symmetric spectra
\beg{ebar0}{\diagram
M_{10}\rto^\gamma & M_{11}\\
M_{00}\uto^\beta\rto_\alpha & M_{01}\uto_\delta
\enddiagram
}
and suppose there exists a map of $R$-modules $s:M_{11}\r M_{10}$ such that
$$\gamma s=Id,$$
$$\beta\vee s:M_{00}\vee M_{11}\r M_{10} \;\text{is an equivalence}.$$
Then 
$$C_2 M\simeq \Sigma M_{01}.$$
\end{lemma}

\Proof
The commutative diagram
\beg{ebar1}{\diagram
M_{00}\vee M_{11}\rto^(.6){\gamma\circ(\beta\vee s)} &M_{11}\\
M_{00}\uto^{\iota_0}\rto_{\alpha} &M_{01}\uto_\delta
\enddiagram
}
maps into \rref{ebar0} by the map $\beta\vee s$ in the upper left corner
and identity elsewhere, and hence has an equivalent $2$-cofiber, since
$\beta\vee s$ is an equivalence. Now since $\gamma s=Id$,
the diagram \rref{ebar1} maps into
\beg{ebar2}{\diagram
M_{00}\rto& 0\\
M_{00}\uto^{Id}\rto_\alpha & M_{01}\uto
\enddiagram
}
with cofiber
$$\diagram
M_{11}\rto^{Id} &M_{11}\\
0\uto\rto & 0\uto,
\enddiagram
$$
so the $2$-cofiber of \rref{ebar2} is equivalent to the $2$-cofiber of \rref{ebar1}.
But \rref{ebar2} in turn maps into
$$\diagram
M_{00}\rto & 0\\
M_{00}\uto^{Id}\rto & 0\uto
\enddiagram
$$
with fiber
\beg{ebar3}{\diagram
0\rto & 0\\
0\uto\rto & M_{01}\uto.
\enddiagram
}
So the $2$-cofiber of \rref{ebar3} is equivalent to the $2$-cofiber of \rref{ebar0}.
\qed

\vspace{3mm}

\begin{lemma}
\label{lbar2}
Consider a diagram $N$ of the form
\beg{ebar20}{\diagram
&N_{101}\rrto^\gamma &&N_{111}\\
N_{001}\urto^\beta\rrto_(.4)\alpha&&N_{011}\urto_\delta &\\
&N_{100}\uuto_(.4)\epsilon\rrto^(.4)\mu && N_{110}\uuto_\pi\\
N_{000}\uuto^\zeta\urto^\eta \rrto_\nu&&N_{010}.\uuto_(.4)\kappa\urto_\lambda&
\enddiagram
}
Assume there exists a map $s:N_{111}\r N_{101}$ such that
$$\gamma s=Id,$$
$$\beta\vee s: N_{001}\vee N_{111}\r N_{101} \;\text{is an equivalence}$$
and assume further that there exists a map 
$$t:N_{100}\r N_{001}$$
such that
$$\beta t=\epsilon,\; \zeta=\eta t.$$
Then 
$$C_3N\simeq \Sigma C_2 M$$
where $M$ is the diagram
\beg{ebar2+}{\diagram
N_{100}\rto^(.4){\mu\Pi\alpha t} &N_{110}\Pi N_{011}\\
N_{000}\uto^\eta\rto_\nu & N_{010}\uto_{\lambda\Pi\kappa}.
\enddiagram
}
\end{lemma}

\Proof
Into \rref{ebar20}, there maps
\beg{ebar22}{\diagram
&N_{001}\vee N_{111}\rrto^{\gamma\circ(\beta\vee s)} &&N_{111}\\
N_{001}\urto^{\iota_1}\rrto_(.4)\alpha&&N_{011}\urto_\delta &\\
&N_{100}\uuto_(.4){\iota_1 t}\rrto^(.4)\mu && N_{110}\uuto_\pi\\
N_{000}\uuto^\zeta\urto^\eta \rrto_\nu&&N_{010}.\uuto_(.4)\kappa\urto_\lambda&
\enddiagram
}
where the map on the $101$-corner is the equivalence
$$\beta\vee s:N_{001}\vee N_{111}\r N_{101}.$$
As in Lemma \ref{lbar1}, \rref{ebar22} maps into
\beg{ebar23}{\diagram
&N_{001}\rrto &&0\\
N_{001}\urto^{Id}\rrto_(.4)\alpha&&N_{011}\urto &\\
&N_{100}\uuto_(.4){ t}\rrto^(.4)\mu && N_{110}\uuto\\
N_{000}\uuto^\zeta\urto^\eta \rrto_\nu&&N_{010}.\uuto_(.4)\kappa\urto_\lambda&
\enddiagram
}
with fiber
$$\diagram
&N_{111}\rrto^{Id} &&N_{111}\\
0\urto\rrto&&0\urto &\\
&0\uuto\rrto && 0\uuto\\
0\uuto\urto \rrto&&0.\uuto\urto&
\enddiagram
$$
By a standard construction in homotopy theory, a diagram of the form \rref{ebar23} can be ``folded'' into the
suspension of the diagram
\beg{ebar24}{\diagram
(N_{100}\Pi N_{001})^\prime\rto^\phi & N_{110}\Pi N_{011}\Pi \widetilde{N_{001}}\\
(N_{000})^\prime\uto^{(\eta\Pi s)^\prime}\rto_\nu & N_{010}\uto_{\lambda\Pi\kappa
\Pi 0}
\enddiagram
}
where $(?)^\prime$ denotes cofibrant replacement and $\widetilde{(?)}$ denotes
fibrant replacement, and $\phi$ is the product of
$$\diagram N_{100}\Pi N_{001}\rto^{\mu p_1} & N_{110}\enddiagram$$
$$\diagram N_{100}\Pi N_{001}\rto^{\alpha p_2} & N_{011}\enddiagram$$
$$\diagram( N_{100}\Pi N_{001})^\prime\rto^(.6){t-Id} & \widetilde{N_{001}}.\enddiagram$$
The diagram \rref{ebar24} commutes up to homotopy, but can be converted into a strict
diagram by standard techniques (for example, by Theorem 1.4 of \cite{em}). (Note:
these complications are, of course, caused by the fact that the canonical
map $(Id\vee 0) \Pi (0\vee Id):A\vee B\r A\Pi B$ is an equivalence but not
an isomorphism in the category of $R$-modules.)

Now into \rref{ebar24}, there maps
$$\diagram(N_{001})^\prime \rto^{Id} &(N_{001})^\prime\\
0\uto \rto & 0\uto
\enddiagram$$
where the upper left corner maps by $0\Pi Id$,
the upper right corner by
$$0\Pi Id,$$
and the upper right corner by
$$0\Pi \alpha \Pi Id$$
(omitting fibrant and cofibrant replacements from the notation). The cofiber
is equivalent to \rref{ebar2+}.\qed

\vspace{3mm}

\noindent
{\bf Proof of Theorem \ref{tinv}:}

As usual, it suffices to prove invariance under R1-moves,
R2-moves and R3-moves.

{\bf Invariance under the R1 move:} After performing an R1-move, 
consider the restrictions of the lax functor $\mathcal{L}_K\mc$ to the 
subcategory (enriched over groupoids) where there is a $0$-resolution
(resp. $1$-resolution) of the new crossing created by the move. The
cobordism from the $0$-resolution to the $1$-resolution will give a
lax natural transformation $\eta$ between these functors. Denote these
restrictions by $\mathcal{L}_K\mc_\epsilon$, $\epsilon =0,1$. Depending on the sign
of the move (which is by definition the sign of the new crossing), 
one of the resolutions will have an extra boundary component
(the $0$-resolution in case of the negative move and the $1$-resolution
in case of the positive move). The new boundary component can be labeled $1$ or
$x$, and this makes this functor $\mathcal{L}_K\mc_\epsilon$
laxly isomorphic to two copies of the functor $\mathcal{L}_K\mc_{1-\epsilon}$.
Further, we can laxly split $\eta$ by choosing this label to be $1$ (in case of the negative
move) and by forgetting the label (in case of the positive move).
In either case, after applying the Elmendorf-Mandell machine, the
cofiber of the realization of $\eta$ becomes isomorphic to the realization
of the other factor of  $\mathcal{L}_K\mc_\epsilon$ resp. its suspension, i.e.
the invariant before the move resp. its suspension, depending
on whether the move was negative or positive.

\vspace{3mm}
{\bf Invariance under the R2 move:} We use the ``Khovanov bracket'' notation of 
Figure 2 of Bar-Natan \cite{dbar}, omitting the suspensions (see Figure 2).

\begin{figure}[!ht]
\includegraphics{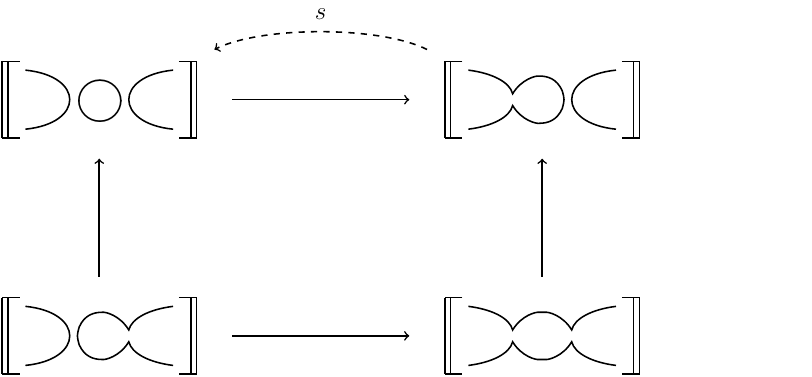}
\caption{}
\end{figure}

We give this picture, however, a modified interpretation: Each bracket denotes a
lax functor $(I^n)^\prime \r\mathcal{S}_2$ corresponding to the indicated partial
resolution of the projection after the R2-move. The arrows in Figure 2
are lax natural transformations.  With the notation
of Lemma \ref{lbar1}, the functor $s$ multiplies objects by the
label $1$ on
the additional connected boundary component. On $1$-morphisms, the functor $s$ tensors a morphism
with $\C$,and $2$-isomorphisms with $Id$. Upon applying the Elmendorf-Mandell machine,
including Theorem 1.4 of \cite{em}, we can obtain a strict functor
$$\diagram
\cdot \rrto && \cdot \lltou_s \\
\cdot\uto\rrto && \cdot \uto
\enddiagram \times I^n\r S-modules$$
which, up to equivalence, has the form
$$\diagram
M\Pi M\rrto_{p_2} && M \lltou_{Id\Pi 0} \\
M\uto^{0\Pi Id}\rrto && ?, \uto
\enddiagram $$
which implies the assumptions of Lemma \ref{lbar1}. Here, $M$ is as in Lemma \ref{lbar1}, and the
$?$ in the lower right corner is the argument of the functor.

\begin{figure}
\includegraphics{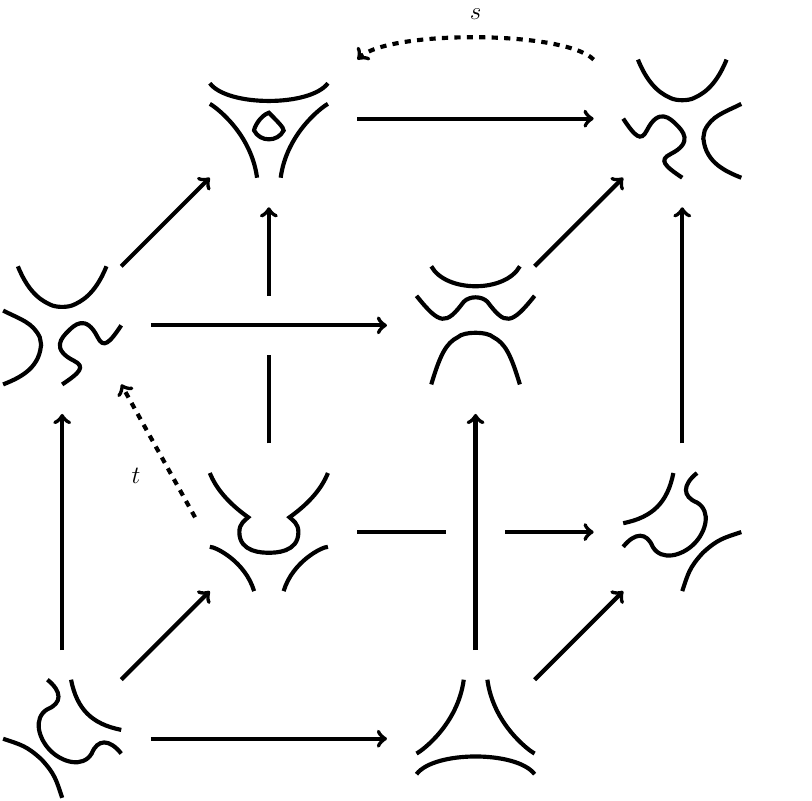}
\caption{}
\end{figure}

\vspace{3mm}
{\bf Invariance under the R3 move:} We follow, again, 
Bar-Natan \cite{dbar}, 
adapting the proof to categories enriched in groupoids. In Figure 3, with the 
notation of Lemma \ref{lbar2},
the constuction of the $s$-map in the $??1$-square is precisely the same
as in the above treatment of the R2-move. Regarding the map $t$, note
that the lax functors $(I^{n})^\prime \r\mathcal{A}_K$ at the $001$ and
$100$ corners are canonically isomorphic (as are the partial resolutions
drawn); let $t$ be the canonical lax isomorphism. From this point on, apply
the Elmendorf-Mandell machine, and use Lemma \ref{lbar2}.
\qed

\vspace{3mm}

\noindent
{\bf Proof of Theorem \ref{t1}:} All that remains to show is that applying
$?\wedge_k H\Z$ to our construction produces an $H\Z$-module which, using
the equivalence \cite{ekmm} Section IV.2, produces a chain complex whose
homology is Khovanov homology. 

To prove this, we note that the strict symmetric ring
spectrum unit 
$$S\r H\Z$$
is realized, on the level of bipermutative categories, by the functor
$$\mathcal{S}_2\r \Z$$
which assigns to a finite set its cardinality. We conclude that applying
$?\wedge H\Z$ to our invariant is realized by taking the Khovanov cube functor
as mentioned in \cite{dbar}, and then applying to it the Elmendorf-Mandell
machinery instead of the totalization described in \cite{dbar}. (Smashing
with $H\Z$ commutes with the Elmendorf-Mandell machine and with the
iterated homotopy cofiber, since it is a left adjoint)

One must, therefore, show that Elmendorf-Mandell machinery \cite{em} to a diagram 
$D$ of
abelian groups (=$\Z$-modules) produce an $H\Z$-module corresponding, under
the machinery of \cite{ekmm} IV.2, to the chain complex obtained as the
homotopy colimit of the diagram $D$ in the category of chain complexes. This follows
from the fact that the equivalence \cite{ekmm} IV.2 commutes, up to equivalence,
with simplicial realization.
\qed


\vspace{10mm}
\begin{thebibliography}{99}

\bibitem{as} M.Atiyah, G.Segal: Twisted K-theory, {\em Ukr. Mat. Visn.}  1  
(2004),  no. 3, 287--330;  translation in  Ukr. Math. Bull.  1  (2004),  no. 3, 291-334

\bibitem{bk} B.Bakalov, A.Kirillov: {\em Lectures on tensor categories and modular functors}, 
University Lecture Series, 21. American Mathematical Society, Providence, RI, 2001


\bibitem{dbar} Dror Bar-Natan: On Khovanov's categorification of the Jones
polynomial, {\em Alg. Geom. Topology} 2 (2002) 337-370

\bibitem{dbar2} Dror Bar-Natan: Khovanov's homology for tangles and cobordisms, {\em Geom. Top.} 9
(2005) 1443-1499

\bibitem{coward} A. Coward: Ordering the Reidemeister moves of a 
classical knot, {\em Alg. Geom. Topol.} 6 (2006) 659-671

\bibitem{ekmm} A.D.Elmendorf, I.Kriz, M.A.Mandell, J.P.May: 
{\em Rings, modules and algebras in stable homotopy theory. With
an Appendix by M.Cole}, Mathematical Surveys and Monographs 47,
American Mathematical Society, Providence, RI, 1997

\bibitem{em} A.D.Elmendorf, M.A.Mandell: 
Rings, modules, and algebras in infinite loop space theory,  
{\em Adv. Math.}  205  (2006),  no. 1, 163-228

\bibitem{fhk} T.M. Fiore, P. Hu, I. Kriz: Laplaza sets, or how to 
select coherence diagrams for pseudo algebras,
{\em Adv. Math.} 218 (2008), no. 6, 1705-1722

\bibitem{kh2} I.Frenkel, M.Khovanov: Canonical bases in tensor products and
graphical calculus for $U_q(sl_2)$, {\em Duke Math. J.}  87  (1997),  no. 3, 409-480



\bibitem{hk} P.Hu, I.Kriz: Conformal field theory and elliptic cohomology,  
{\em Adv. Math.}  189  (2004),  no. 2, 325-412

\bibitem{huang} Y.Z.Huang: Vertex operator algebras, the Verlinde conjecture, 
and modular tensor categories.  
{\em Proc. Natl. Acad. Sci. USA}  102  (2005),  no. 15, 5352-5356

\bibitem{huang1} Y.Z.Huang: Vertex operator algebras and the Verlinde conjecture.  
{\em Commun. Contemp. Math.}  10  (2008),  no. 1, 103-154

\bibitem{jardine} J.F.Jardine: Cubical Homotopy Theory, a beginning, {\em Isaac Newton Institute preprints}, 2002


\bibitem{khovanov} M.Khovanov: A categorification of the Jones polynomial,
{\em Duke Math. J.} 101 (2000) 359-426

\bibitem{kh3} M.Khovanov: $sl(3)$ link homology,  {\em
Algebr. Geom. Topol.}  4  (2004), 1045-1081

\bibitem{kk} D.Kriz, I.Kriz: Baldwin-Ozsv\'{a}th-Szab\'{o} cohomology is
a link invariant, arXiv:1109.0064

\bibitem{spin} I.Kriz: On spin and modularity in conformal field theory,  
{\em Ann. Sci. Ecole Norm. Sup.} (4) 36 (2003), no. 1, 57-112

\bibitem{lls} T.Lawson, R.Lipshitz, S.Sarkar: 
Khovanov homotopy type, Burnside category, and products , arXiv:  arXiv:1505.00512

\bibitem{ls} R.Lipshitz, S.Sarkar: A Khovanov homotopy type or two, arXiv:1112.3932




\bibitem{klai} I.Kriz, L.Lai: On the definition and K-theory realization of a modular functor,  arXiv:1310.5174

\bibitem{lurie} J.Lurie: Expository article on topological field theories, preprint, 2009




\bibitem{ms} G.Moore, N.Seiberg: Classical and quantum conformal field theory,  
{\em Comm. Math. Phys.}  123  (1989),  no. 2, 177-254

\bibitem{os} P.Ozsv\'{a}th, Z.Szab\'{o}: Holomorphic disks and knot invariants.  
{\em Adv. Math.} 186  (2004),  no. 1, 58-116

\bibitem{os1} P.Ozsv\'{a}th, Z.Szab\'{o}: On the Heegaard Floer homology 
of branched double-covers, {\em Adv. Math.}  194  (2005),  no. 1, 1-33



\bibitem{scft} G.Segal: The definition of conformal  field theory, 
{\em  Topology, geometry and quantum field theory},  421-577, 
London Math. Soc. Lecture Note Ser., 308, Cambridge Univ. Press, Cambridge, 2004

\bibitem{verlinde} E. Verlinde: Fusion rules and modular transformations in 2D 
conformal field theory, {\em Nuclear Phys. B}  300  (1988),  no. 3, 360-376

\bibitem{cat} B.Webster: Knot invariantes and higher representation theory I,II,
arXiv: 1001.2020, 1005.4559


\end{thebibliography}
\end{document}